\documentclass[12pt]{amsart}
\usepackage{mathrsfs}
\usepackage{amsfonts}
\usepackage{amsthm}
\usepackage{mathrsfs}
\newtheorem{theorem}{Theorem}
\newtheorem{lemma}{Lemma}
\newtheorem{corollary}{Corollary}

\newtheorem{proposition}{Proposition}

\newcommand{\C}{\mathbb{C}}
\newcommand{\D}{\Omega}

\begin{document}

\title[Non-compactness of the $\overline{\partial}$-Neumann
operator]{Analytic discs, plurisubharmonic hulls, and non-compactness
of the $\overline{\partial}$-Neumann operator}

\author{S\"onmez \c sahuto\u glu}
\author{Emil J. Straube}
\email[S\"onmez \c Sahuto\u glu]{sahutogl@math.tamu.edu}
\email[Emil J. Straube]{straube@math.tamu.edu}  

\address{Department of Mathematics\\ Texas A\&M University\\
              College Station, TX 77843-3368}

\thanks{Research supported in part by NSF grant number DMS-0100517.} 
 

\begin{abstract}
We show that a complex manifold $M$ in the boundary of a smooth
bounded pseudoconvex domain $\D$ in $\C^{n}$ is an obstruction to
compactness of the $\overline{\partial}$-Neumann operator on $\D$,
provided that at some point of $M$, the Levi form of $b\D$ has the
maximal possible rank $\, n-1-dim(M)$ (i.e. the boundary is strictly
pseudoconvex in the directions transverse to $M$). In particular, an
analytic disc is an obstruction, provided that at some point of the
disc, the Levi form has only one zero eigenvalue (i.e. the eigenvalue
zero has multiplicity one). We also show that a boundary point where
the Levi form has only one zero eigenvalue can be picked up by the
plurisubharmonic hull of a set only via an analytic disc in the
boundary.
\end{abstract}

\maketitle 

\section{Introduction and results}

Let $\D$ be a smooth bounded pseudoconvex domain in $\C^n.$ The
$\overline{\partial}$-Neumann operator $N$ is the inverse of the
complex Laplacian $\Box=\overline{\partial}\overline{\partial}^*+
\overline{\partial}^*\overline{\partial}$ on square integrable
$(0,1)$-forms on $\D$. Its regularity theory plays an important role
both in partial differential equations and in several complex
variables; detailed information can be found in \cite{BS, CS, FK}. In
particular, the question of when $N$ is, or is not, compact is of
interest in various contexts, among them global regularity (\cite{KN,
BS}), Toeplitz operators (\cite{FS01} and its references),
semiclassical analysis of Schr\"{o}dinger operators (\cite{FS02,
CF}), and kernels solving $\overline{\partial}$ (\cite{HL}). An
important motivation for
studying compactness has emerged over the last ten years, as it
became increasingly clear how subtle the (a priori more important)
property of global regularity is (see the surveys \cite{BS, C}):
compactness is more robust (it localizes, for example) and thus 
should be more readily amenable to a characterization in terms of
properties of the boundary. 

There is a fairly general potential theoretic condition on the
boundary of a domain which implies that the
$\overline{\partial}$-Neumann operator is compact. This condition was
introduced, and shown to imply compactness, by Catlin(\cite{C84})
under the name of property(P). A detailed study of this property,
from the point of view of the Choquet theory of the cone of
plurisubharmonic functions, can be found in \cite{S2}. In particular,
the boundary of a smooth bounded pseudoconvex domain satisfies
property(P) if and only if continuous functions on the boundary can
be approximated uniformly (on the boundary) by functions continuous
on the closure and plurisubharmonic in the interior (\cite{S2},
th\'{e}or\`{e}me 2.1). For a discussion of property(P) in the spirit of Oka's lemma and related results, see \cite{Ha04}. McNeal showed that a relaxed version of property(P) still implies compactness (\cite{McN}). An interesting sufficient condition for compactness, related to both property(P) and the relaxed version in \cite{McN}, was given in \cite{T91}. Recently, the second author introduced a more geometric approach to compactness of the $\overline{\partial}$-Neumann operator on domains in $\C^{2}$ (\cite{S04}). 

The most blatant violation of propertiy(P) as well as
of the sufficient conditions in \cite{T91}, \cite{McN}, and \cite{S04} is an analytic disc in the boundary. Such discs also cause failure of hypoellipticity for $\overline{\partial}$ \cite{C81, DP}). The
question whether they also cause failure of compactness is thus a
natural one. That the answer is yes in the case
of domains in $\C^{2}$ is an unpublished result of Catlin from the
early eighties; a proof for the case of domains with Lipschitz
boundary is in \cite{FS01}. (Note that the converse is false: there
are obstructions to compactness of the $\overline{\partial}$-Neumann
operator considerably more subtle than discs in the boundary, see
\cite{M97, FS01}.) The situation is completely understood in locally
convexifiable domains: $N$ is compact (if and) only if the boundary
contains no analytic disc (\cite{FS98, FS01}). It is folklore that
the methods that work in $\C^{2}$ also show that, in any dimension
$n$, if the $\overline{\partial}$-Neumann operator is compact, then
the boundary cannot contain an $(n-1)$-dimensional complex manifold.
However, whether a disc is necessarily an obstruction to compactness
is open. We show in this paper that this is the case when the disc
contains a point at which the boundary is strictly pseudoconvex in
the directions transverse to the disc. More generally, one can trade
one strictly positive eigenvalue of the Levi form for an increase of
one in the dimension of the complex manifold that is to obstruct
compactness.

Although we are mainly interested in compactness of the
$\overline{\partial}$-Neumann operator $N$, our arguments proceed via
existence of a compact solution operator to $\overline{\partial}$ on
$(0,1)$-forms. A solution operator $T$ for $\overline{\partial}$ on
$(0,1)$-forms is a (continuous) operator $T : L^{2}_{(0,1)}(\Omega)
\cap ker(\overline{\partial}) \rightarrow L^{2}(\Omega)$ such that
$\overline{\partial}T(\alpha) = \alpha$ for all $\alpha \in
L^{2}_{(0,1)}(\Omega) \cap ker(\overline{\partial})$. It is well
known that if $N$ is compact, then the canonical solution operator
$\overline{\partial}^{*}N$ is compact (considerably more can be said,
see Lemma 1.1 in \cite{FS01}). Note that since the projection of
$L^{2}(\Omega)$ onto the orthogonal complement of the holomorphic
functions preserves compactness, saying that there exists a compact
solution operator for $\overline{\partial}$ is the same as saying
that the canonical solution operator is compact.

\begin{theorem}\label{thm6}
Let $\D$ be a smooth bounded pseudoconvex domain in $\C^n, n\geq 2.$
Let $P \in b\D$ and assume that the Levi form of $b\D$ at $P$ has the
eigenvalue zero with multiplicity at most $k, 1 \leq k \leq n-1$. If
there exists a compact solution operator for $\overline{\partial}$ on
$(0,1)$-forms (in particular, if the $\overline{\partial}$-Neumann
operator on $\D$ is compact), then $b\D$ does not contain a
$k$-dimensional complex manifold through $P$.
\end{theorem}

In particular, if the Levi form of the boundary has at most one
degenerate eigenvalue at all boundary points, then compactness of the
$\overline{\partial}$-Neumann operator implies that there are no
discs in the boundary. A special case of this occurs
imlicitly in \cite{K}, Theorem 1.1, for domains fibered over a
Reinhardt domain in $\C^{2}$. We will give the proof of Theorem
\ref{thm6} in section 2, along with the proof of the next corollary.
Note that when the set of \emph{Levi flat} boundary points has
nonempty (relative) interior, then it is foliated by
$(n-1)$-dimensional complex manifolds. By the folklore result
mentioned above, this is incompatible with compactness of the
$\overline{\partial}$-Neumann operator. It is an easy consequence of
Theorem \ref{thm6} that this result holds for the (in general much
bigger) set of \emph{weakly pseudoconvex} points.

\begin{corollary}\label{col1}
Let $\D$ be a smooth bounded pseudoconvex domain in $\C^n,n\geq 2
$. If there exists a compact solution operator for
$\overline{\partial}$ on $(0,1)$-forms (in particular, if the
$\overline{\partial}$-Neumann operator on $\D$ is compact), then the
set of weakly pseudoconvex points in $b\D$ has empty (relative)
interior.
\end{corollary}

\emph{Remark 1.} If $b\Omega$ is not strictly pseudoconvex in the
directions transverse to the complex submanifold at $P$, then the
boundary is `more flat'. In general, one would expect this situation
to be even more favorable to noncompactness of the
$\overline{\partial}$-Neumann operator. However, our present methods
do not give this.
\medskip

Denote by $P(\overline{\D})$ the set of
plurisubharmonic functions on $\D$ that are continuous up to the
boundary of $\D.$  Define the plurisubharmonic hull, $\widehat{K},$
of a compact set $K\subset \overline{\D}$ as follows: 
$$\widehat{K}=\left\{ z\in \overline{\D}:
f(z)\leq \sup_{w\in K}f(w) \,\,{\rm for\,\,all}\,\,f\in
P(\overline{\D}) \right\}$$
The same hull results when taken with respect to any of the following
sets of functions: plurisubharmonic functions smooth up to the
boundary, holomorphic functions smooth up to the boundary, or
holomorphic functions continuous up to the boundary (see
for example \cite{C78,C80,HS}). The intersection of $\widehat{K}$ with
the boundary is determined by $K \cap b\D$: $\widehat{K} \cap b\D =
(\widehat{K \cap b\D}) \cap b\D$ (\cite{C80}, Theorem 3.1.7,
\cite{HS} Proposition 2). Thus a compact subset $K$ of
$\overline{\D}$ which picks up hull in the boundary, that is, where  
$\widehat{K} \cap b\D \neq K \cap b\D$, represents an obvious
violation of property(P) (the situation with respect to the conditions
in \cite{T91}, \cite{McN}, and \cite{S04} seems not to be understood at
present) as well as an
obstruction to hypoellipticity of $\overline{\partial}$ (\cite{C81},
Theorem 3). Such sets therefore present themselves as interesting
candidates for obstructions to compactness in the
$\overline{\partial}$-Neumann problem. In $\C^{2}$, this leads to
nothing new: a compact $K$ can pick up hull in the boundary only via
an analytic disc contained in the boundary (\cite{C81}, Theorem 4).
In higher dimensions, however, sets can pick up hull in the boundary
in less obvious ways, that is, even when the boundary contains no
analytic disc (\cite{C81}, Theorem 5). Thus, in view of the case
$k=1$ of Theorem \ref{thm6}, the question arises how a point in the
boundary can be picked up by the plurisubharmonic hull if the Levi
form at the point has only one degenerate eigenvalue. Perhaps not
surprisingly, the situation is the same as in $\C^{2}$: there has to
be an analytic disc in the boundary through the point. 

\begin{proposition}\label{hull}
Let $\D$ be a smooth bounded pseudoconvex domain in $\C^n,n\geq 2 .$
Let $P$ be a boundary point where the Levi form of $b\D$ has at most
one degenerate eigenvalue. Then there is a compact subset $K$ of
$\overline{\D}$ with $P \in \widehat{K} \setminus K$ if and only if
$b\D$ contains an analytic disc through $P$.
\end{proposition}

Proposition \ref{hull} can be shown along the lines of the proof for
the $\C^{2}$ case given in \cite{C78}, but some additional work is
needed. We will give the details in section 3.

\bigskip

\emph{Remark 2.} Sibony has characterized the domains for which
$\widehat{K} \cap b\D = K \cap b\D$ for all compact subsets $K$ of
$\overline{\D}$ by hypoellipticity of $\overline{\partial}$ in
\emph{local} $L^{2}$-spaces: for every
$\overline{\partial}$-closed (0,1)-form $\alpha$ with locally square
integrable coefficients in $\D$ there should exist a locally square
integrable function $u$ on $\D$ solving $\overline{\partial}u =
\alpha$ and such that the singular support of $u$ in $\overline{\D}$
is equal to the singular support of $\alpha$ (\cite{S1}). Combining
this with Theorem \ref{thm6} and Proposition \ref{hull} above shows
that on a
smooth bounded pseudoconvex domain, if the Levi form has at most one
degenerate eigenvalue at each boundary point, compactness of the
$\overline{\partial}$-Neumann operator implies hypoellipticity of
$\overline{\partial}$ in the space of locally square integrable
forms. It would be very interesting to have a direct proof of this
implication.

\section{Proofs of Theorem \ref{thm6} and Corollary \ref{col1}} 
We will need to control the boundary geometry near a point of a
complex manifold in the boundary. This is accomplished in the
following lemma.

\begin{lemma} \label{lem1}
Let $\D$ be a smooth bounded pseudoconvex domain in $\C^n$, $M$ a
complex 
manifold of dimension $k$ in $b\D,$ and $P\in M.$ Then there is a ball
$B$ 
centered at $P,$ a biholomorphic map $G:B\to G(B)$ such that
\begin{itemize}
\item[(i)] $G(P)=0$
\item[(ii)] $G(M\cap B)=\{w\in G(B)|w_{k+1}= \cdots =w_n=0\}$
\item[(iii)] the real normal to $G(b\D \cap B)$ at points of $G(M\cap
B)$ is 
given by the $Re(w_n)$-axis. 
\end{itemize} 
\end{lemma} 
In other words, there is a local holomorphic change of coordinates so
that in the new coordinates, $M$ is affine and the (real) unit normal
to the boundary is constant along $M$. 

To prove the lemma, note first that we can change coordinates near $P$
so that 
$M$ is locally given by $\{z|z_{k+1}= \cdots =z_n=0 \},$ and $P=0.$ In
these coordinates, 
let $\rho$ be a defining function for $b\D$ near 0. We may assume that
$\partial \rho/\partial z_n \neq 0$. There exists a
real-valued $C^{\infty}$
function $h$ in a neighborhood of 0 such that the normal 
$e^h\left(0,\cdots,0, \frac{\partial \rho}{\partial \bar
z_{k+1}},\cdots, 
\frac{\partial \rho}{\partial \bar z_n} \right)$
is conjugate holomorphic on $M$ (each component is conjugate
holomorphic). When $M$ is a 
disc, this is the conclusion of Lemma 1 in \cite{BF81}. One can adapt
these methods to 
cover the case where $dim M>1.$ Alternatively, the statement is a
special case of 
the theorem in \cite{SS02}, specifically the equivalence of $(ii)$ and
$(iv).$ Note 
that the form $\alpha$ appearing in $(iv)$ of that theorem is real,
and its restriction 
to $M$ is closed by the lemma in section 2 of \cite{BS93}. Therefore,
there is a real valued 
function $h$ such that  $d(-h)=\alpha$ on $M$ (near 0). The proof of
$(ii)\Leftrightarrow (iv)$
in \cite{SS02} shows that any real-valued $C^{\infty}$ extension of
$h$ to a full neighborhood 
of 0 will do.

We now define a biholomorphic change of coordinates (near 0) by 
$ \widehat{G}(z_1,\cdots,z_n)=(z_1,\cdots,z_n) +S(z)$, where 
$$S(z)= \left( 0,\cdots,0, 
\sum_{j=k+1}^n z_je^{h(z_1,\cdots,z_k,0,\cdots,0)}\frac{\partial
\rho}{\partial z_j}(z_1,\cdots,z_k,0,\cdots,0) \right).$$
Then $\widehat{G}(M)\subset \{z_{k+1}= \cdots  z_n=0\}$, and at points
of $M$, the complex derivative 
of $\widehat{G}$ maps the complex tangent space to $b\D$ onto the
complex hypersurface 
$\{z_n=0\}$. Consequently, in these new coordinates ( which we again
denote by 
$(z_1,\cdots,z_n)$), $M$ is of the same form as before, but the
complex tangent space to 
$b\D$  is constant, namely it is the hyperplane $\{z_n=0\}.$ Consider
now the real unit 
normal to $b\D .$ Its restriction to $M$ is of the form
$(0,\cdots,0,e^{i\theta}).$  
Using Lemma 1 from \cite{BF81} once more shows that the function
$\theta$ is harmonic 
on each disc in $M$, i.e. it is  pluriharmonic on $M.$ Denote by $h_1$
a pluriharmonic 
conjugate. The final coordinate change 
$$(z_1,\cdots,z_n)\to \left(
z_1,\cdots,z_{n-1},z{_n}e^{-h_1(z_1,\cdots,z_k,0,\cdots,0)
-i\theta(z_1,\cdots,z_k,0,\cdots,0)} \right)$$ rotates the real
normal so that the unit normal becomes constant on $M$. Combining the
three local biholomorphic coordinate changes gives the map $G$ with
the properties required in Lemma \ref{lem1}.  
\medskip

The proof of Theorem \ref{thm6} will be indirect; we will show that if
the boundary contains a $k$-dimensional complex manifold through $P,$
then the $\overline{\partial}$-Neumann operator on $\D$ is not
compact. The source of noncompactness is contained in the following
lemma. Denote by $A(\D)$ the Bergman space of a domain $\D$ in
$\mathbb{C}^{n}$, that is, the closed subspace of $L^2(\D )$
consisting of holomorphic functions.

\begin{lemma} \label{lem2}
Let $\D$ be a bounded pseudoconvex domain in $\C^n$, smooth near the
strictly pseudoconvex boundary point $P$. Assume the pseudoconvex
domain $\D_1$ is contained in $\D$, shares the boundary point $P$,
and is smooth near $P$. Then the restriction operator from $A(\D)$ to
$A(\D_1)$ is not compact.
\end {lemma}

We prove the lemma by constructing a sequence of norm $1$ functions in
$A(\D)$ which has no convergent subsequence in $A(\D_1)$. Let
$\{P_j\}_{j=1}^{\infty}$ be a sequence of points on the common
interior normal to $b\D$ and $b\D_1$ at $P$ such that $\lim_{j
\rightarrow \infty}P_j = P$. Set
$f_j(z)=K_{\D}(z,P_j)/K_{\D}(P_j,P_j)^{1/2}$, where $K_{\D}$ denotes
the Bergman kernel of $\D$. The normalization is so that $\parallel
f_j\parallel _{L^{2}(\D)} = 1$ for all $j$. For $z \in \D$ fixed, the
function $K_{\D}(z,\cdot)$ is smooth up to the boundary near $P$ (in
fact, the subelliptic estimates for the $\overline{\partial}$-Neumann
problem near $P$ have considerably stronger consequences for the
kernel function \cite{B87}, \cite{B86}), and $K_{\D}(P_j,P_j)
\rightarrow \infty$ as $j \rightarrow \infty$ (see e.g. \cite{H65},
Theorem 3.5.1). Consequently, $f_j(z) \rightarrow 0$ for all $z \in
\D$. On the other hand,
\[\parallel f_j\parallel _{L^{2}(\D_1)}^{2}=\frac{\parallel
K_{\D}(\cdot,P_j)\parallel _{L^{2}(\D_1)}^{2}}{K_{\D}(P_j,P_j)}
\geq \frac{K_{\D}(P_j,P_j)}{K_{\D_1}(P_j,P_j)} \geq C > 0\ . \]
The first inequality follows by applying the reproducing property of
$K_{\D_1}$ to the function $K_{\D}(\cdot,P_j)$ viewed as an element
of $A(\D_1)$. (This argument comes verbatim from \cite{FS98},
p.~637.) Because both $\D$ and $\D_1$ are strictly pseudoconvex at
$P$, the kernel asymptotics at $P$ are the same (\cite{H65}, Theorem
3.5.1), and the second inequality follows. We conclude that
$\{f_j\}_{j=1}^{\infty}$ has no subsequence that converges in
$A(\D_1)$, and the proof of Lemma \ref{lem2} is complete.
\medskip

To prove Theorem \ref{thm6}, we use the simple geometry of $b\Omega$
in the new coordinates (near $P$) to construct a sequence of
$\overline{\partial}$-closed $(0,1)$-forms that will lead to
a contradiction when combined with the existence of a compact
solution operator for $\overline{\partial}$ on $\Omega$. The argument
follows very closely \cite{FS98},
section 4, and \cite{FS01}, proof of Proposition 4.1. In turn, those
arguments draw substantially on ideas from \cite{C81} and \cite{DP}.
We keep the notation from Lemma \ref{lem1}. Denote $G(B \cap \D)$ by
$\widetilde{\D}$. There exist open neighborhoods $V_2\subset\subset
V_1\subset\subset 
V_0\subset\subset G(B)$ of $P=0$ in $\C^n$ and a small enough ball
$\widetilde B$ in $\C^{n-k}$ such that 
$M_0\times \widetilde{B} \subset \widetilde{\D}$, where we set
$M_j=M\cap V_j, j=0,1,2$. Here we have used that the real normal to
the boundary is constant along $M$. Let $S_0$ be 
the slice of $V_0\cap \widetilde{\D}$ in $\C^{n-k}$ that is
perpendicular to $M$. Because the rank of the Levi form is preserved
under biholomorphisms, $S_0$ is strictly pseudoconvex at $0$.
Lemma \ref{lem2} shows that there exists a sequence of holomorphic
functions $\{ f_j\}_{j=1}^{\infty}$ that is bounded in $A(S_0)$
and has no convergent subsequence in $A(\widetilde{B}).$  Using the
Ohsawa-Takegoshi extension theorem \cite{OT} one can extend $f_j$ from
$S_0$ to $\widetilde{\D}$  for each $j$ to get a bounded sequence of
holomorphic functions $\{F_j\}_{j=1}^{\infty}$ in 
$A(\widetilde{\D})$. Define
$\beta_k(z)=F_k(z)\chi (z)$
where $ \chi : V_0 \to [0,1]$ is a smooth function such that
 \begin{displaymath}
\chi (z) = \left\{ \begin{array}{ll}
1 & \textrm{if}\quad z\in  M_2\times \widetilde{B} \\
0 & \!\! \textrm{ if} \quad z\not\in V_1\\
\end{array} \right.
\end{displaymath}
Let $\alpha_k(z)=\overline{\partial} (\chi
(z)F_k(z))=F_{k}(z)\overline{\partial}\chi(z)$.
$\{\alpha_k\}_{k=1}^{\infty}$ is a bounded sequence in
$L_{(0,1)}^2(\widetilde{\D})$. We pull
$\alpha_{k}$ back via $G$ and view the resulting form
$\widetilde{\alpha_{k}}$ as a $\overline{\partial}$-closed form in
$\Omega$ (by extending by $0$ outside $B$). The sequence
$\{\widetilde{\alpha_{k}}\}_{k=1}^{\infty}$ is bounded in
$L^{2}_{(0,1)}(\Omega)$. By assumption, there is a compact solution
operator for
$\overline{\partial}$ on $\Omega$. Applying this solution operator to
the sequence $\{\widetilde{\alpha_{k}}\}_{k=1}^{\infty}$ and
passing to a subsequence if necessary yields a sequence
$\{\widetilde{g_{k}}\}_{k=1}^{\infty}$ which converges in
$L^2(\Omega)$ such that for all $k$, we have
$\overline{\partial}\widetilde{g_{k}} = \widetilde{\alpha_{k}}$.
Restricting to $B \cap \Omega$ and pushing forward to
$\widetilde{\D}$, gives a
sequence $\{g_{k}\}_{k=1}^{\infty}$ that converges in
$L^2(\widetilde{\D})$ and
such that for all $k$, $\overline{\partial}g_{k} = \alpha_{k}$.
Finally, setting $h_k(z)=g_k(z)-\chi (z)F_k(z)$
gives a sequence of  holomorphic functions on $\widetilde{\D}$. Notice
that $h_k=g_k$ on $\widetilde{\D} \setminus V_1$, so
$\{h_{k}\}_{k=1}^{\infty}$ converges in $L^2(\widetilde{\D} \setminus
\overline{V_1})$. The $L^2$-norm of a holomorphic function on a
spherical shell dominates the $L^2$-norm of the function on the whole
ball. Combining this observation with Fubini's theorem implies that
$\{h_k \}_{k=1}^{\infty}$ converges in
$L^2(M_2 \times \widetilde{B})$. Hence so does
$\{F_k\}_{k=1}^{\infty}$ (since $\{g_{k}\}_{k=1}^{\infty}$ also
converges there and $\chi \equiv 1$ on $M_{2} \times \widetilde{B}$).
Using the submean value
property on balls in the variables $w_{1}, w_{2}, \cdots , w_{k}$
while $(w_{k+1}, \cdots , w_{n}) \in \widetilde{B}$ stays fixed and
then again Fubini's theorem we get that $\{f_k\}_{k=1}^{\infty}$
converges in $A(\widetilde{B})$. This is the contradiction we seek.
The proof of Theorem \ref{thm6} is therefore complete.

\medskip

We next prove Corollary \ref{col1}. Let $V$ be a nonempty open subset
of $b\Omega$ contained in the set of weakly pseudoconvex points.
Define $m$ to be the maximum rank of the Levi form on $V$ and let $P
\in V$ be a point where the Levi form has rank $m$ (such a point
exists, since the rank assumes only finitely many values on $V$).
Then near $P$, the rank is at least $m$, hence equal to $m$.
Therefore, $b\Omega$ is foliated, near $P$, by complex manifolds of
dimension $n-1-m$ (see for example \cite{F}). Theorem \ref{thm6}
(for $k=n-1-m$) implies that the $\overline{\partial}$-Neumann
operator on $\Omega$
is not compact. This contradicts the assumption in Corollary
\ref{col1}.

\section{Proof of Proposition \ref{hull}} 

We denote the Levi form of $b\D$ by $L(Y,\overline{Z})$ for complex
tangential vector fields $Y,Z$ of type $(1,0)$. We similarly use the
notation $L_{f}(Y,\overline{Z})$ for the full complex Hessian of a
function $f$. Let $J$ denote the complex structure map of (the tangent
bundle of) $\C^n$, $F_A^t(P)$
the flow generated by the (real) vector field $A$, and
$A^{\theta}=\cos (\theta)A+\sin (\theta )J(A)$. The following lemma
is proved in \cite{C78} for domains in $\C^{2}$.
 
\begin{lemma}\label{flow}
Let $\D$ be a smooth bounded pseudoconvex domain in $\C^n, n\geq 2$,
and $P\in b\D$. Assume that the Levi form of $b\D$ has one
degenerate eigenvalue at $P$. If there exist a neighborhood $U$ of $P$
in $b\D$, a smooth complex tangential (real) vector field $A$ defined
on $U,$ and $t_0>0$  such that $L(A-iJ(A))=0$ at all points of the
real two-dimensional manifold $F_{A^{\theta }}^t(P)$, $0\leq t \leq
t_0$, $0\leq \theta \leq 2\pi ,$ then this manifold is an analytic
disc, possibly after shrinking $t_0$.
\end{lemma}

We want to show that the real two-dimensional manifold
$M=\{F_{A^{\theta}}^t(p):0\leq t \leq \widehat{t_0},0\leq \theta \leq
2\pi\}$ is actually a one dimensional complex manifold for some
$\widehat{t_0} > 0$. We follow \cite{C78}
and indicate the necessary modifications. It suffices to show that
the tangent space of $M$ is spanned by $A$ and $J(A)$. This means
that it is invariant under $J$, and therefore $M$ is indeed a complex
manifold. If we view $\theta$ as a parameter in the initial value
problem that determines the flow generated by $A^{\theta}$ and set $f(t,\theta):=F_{A^{\theta}}^t(p)$, we see
that $\frac{\partial f}{\partial \theta}$ is a solution of the
following initial value problem:
\begin{equation}\label{eq1} 
\begin{cases}
(\frac{\partial f}{\partial \theta})^{\prime}=\frac{\partial
A^{\theta}}{\partial x}\frac{\partial f}{\partial \theta} +
\frac{\partial A^{\theta}}{\partial \theta}\\
	\frac{\partial f}{\partial \theta}(0)=0
\end{cases}
\end{equation}
where $\frac{\partial A^{\theta}}{\partial x}$ is the Jacobian of
$A^{\theta}$. We will find a solution to
(\ref{eq1}) in the form $a(t)A+b(t)J(A)$. Then $\frac{\partial
f}{\partial \theta}$, the unique solution of (\ref{eq1}), is of this
form. Therefore, both $\frac{\partial f}{\partial t}$ and
$\frac{\partial f}{\partial \theta}$ are linear combinations of $A$
and $J(A)$. This will complete the proof of Lemma \ref{flow}.

If we substitute $a(t)A+b(t)J(A)$ for $\frac{\partial
f}{\partial \theta}$ into (\ref{eq1}), we get:
\begin{eqnarray}\label{eq2}
 0&=&	a'(t)A+(a(t)\sin(\theta)-b(t)\cos (\theta))[J(A),A]\\
 &&+b'(t)J(A)+\sin (\theta)A-\cos(\theta)J(A)\nonumber \\
	0&=&a(0)=b(0) \nonumber
\end{eqnarray}
Suppose that we can show that $[J(A),A]$ is a linear combination of
$A$ and $J(A)$ at points of $M$. Note that the coefficients are
automatically smooth. Then after collecting terms containing $A$ and
$J(A)$, respectively, \eqref{eq2} becomes an initial value problem
for $(a(t),b(t))$ that has a (unique) solution, and we are done. It is
clear that $[J(A),A]$ is complex tangential (at points of $M$),
because $A-iJ(A)$ is a Levi null direction. In $\C^{2}$, this means
that it is a (real) linear combination of $A$ and $J(A)$. The
argument so far comes entirely from \cite{C78}. In higher dimensions,
some additional work is needed. 

The hypothesis that the Levi form has one degenerate eigenvalue at $P$
(hence at most one near $P$) implies that it is diagonalizable near
$P$ (see e.g. \cite{M}, Lemma 2.1).
Let $X_1,X_2,\cdots,X_{n-1}$ be smooth complex tangential vector
fields of type $(1,0)$ that diagonalize the Levi form near $P$. We
may assume that $X_1= A-iJ(A)$ on $M$ (near $P$), since $A-iJ(A)$ is
an eigenvector associated with the eigenvalue $0$ at points of $M$.
Then the flows $F^{t}_{B^{\theta}}(P)$ associated with the vector
fields $B^{\theta}$, where $B=\mathrm{Re}X_1$, also generate $M$, and
we may therefore assume that $A-iJ(A) = X_1$ in a full neighborhood
of $P$ in the boundary. (In the proof of Proposition \ref{hull} below,
we will actually be in this situation from the outset.) Thus the
commutator we are interested in is $\left[X_1,\overline{X_1}\right] =
2i[B,J(B)]$. If $X_1$
were a Levi null field in a \emph{full} neighborhood of $P$, then
$[X_1, \overline{X_1}]$ would likewise be (\cite{F}). However, we
only know that $L(X_1, \overline{X_1}) =0$ on $M$. What we can assert
is that $\left[X_1, \overline{X_1}\right] = Y-\overline{Y}+\varphi
(L_n-\overline{L_n})$, where $Y$ is a smooth complex tangential field
of type $(1,0)$, $L_n$ is the complex normal to the boundary, and
$\varphi$ is a smooth, nonnegative function that vanishes on $M$. The
nonnegativity of $\varphi$ is a consequence of the pseudoconvexity of
$\Omega$. What we need is that on $M$, $Y$ is a multiple of $X_1$.
The Jacobi identity for $X_1$, $\overline{X_1}$, and $X_k$ gives 
\begin{equation}
\left[Y-\overline{Y}+\varphi (L_n-\overline{L_n}),X_k\right] +
\left[[\overline{X}_1,X_k],X_1\right]+\left[[X_k,X_1],\overline{X}_1
\right]=0,k=1,\cdots,n-1; \nonumber
\end{equation} 
we have replaced $\left[X_1, \overline{X_1}\right]$ in the first term
by  $Y-\overline{Y}+\varphi (L_n-\overline{L_n})$. The second and
third commutators are complex tangential at points of $M$ if $k \geq
2$. To see this, note that $X_1$ is in the nullspace of the Levi form
at points of $M$ and that both the commutators $\left[X_k,X_1\right]$
and $\left[\overline{X_1},X_k\right]$ are complex tangential in a
\emph{full} neighborhood of $P$ in $b\Omega$ (the latter because
$X_1, \cdots ,X_{n-1}$ diagonalize the Levi form near $P$).
Consequently, the first commutator is complex tangential at points of
$M$ as well. In addition, because $\varphi$ and its first order
derivatives vanish at points of $M$, 
$\left[\varphi (L_n-\overline{L_n}),X_k\right]$ is zero at points of
$M$. Therefore, $\left[Y-\overline{Y},X_k\right]$ is complex
tangential at points of $M$, and hence so is
$\left[\overline{Y},X_k\right]$. It follows that at points of $M$,
$Y$ is in the nullspace of the Levi form and thus is a multiple of
$X_1$. This completes the proof of Lemma \ref{flow}. We remark that a
similar use of the Jacobi identity occurs in \cite{F}.
\smallskip
\enlargethispage{20pts}

Once Lemma \ref{flow} is in hand, the proof of Proposition \ref{hull}
can be completed as in \cite{C78}, with only small modifications.
Assume there is a compact subset $K$ of $\overline{\D}$ and a
boundary point $P \in \widehat{K} \setminus K$ and there is no
analytic disc in the boundary through $P$. We may assume that zero is
an eigenvalue (with multiplicity one) of the Levi form at $P$. Take
$V$ to be a neighborhood of $P$ small enough so that $V \cap K =
\emptyset$ and so that there are complex tangential fields $X_1,
\cdots , X_{n-1}$ of type $(1,0)$ which diagonalize the Levi form in
a neighorhood of $\overline{b\D \cap V}$, with
$L(X_1,\overline{X_1})(P) = 0$. By Lemma \ref{flow}, there exist
$\theta$ and $t_0 >0$ such that $F_{A^{\theta}}^{t_0}(P) \in V$ and
$F_{A^{\theta}}^{t_0}(P)$ is a strongly pseudoconvex point, where
$A=X_1+\overline{X_1}$. Near $P$, choose a boundary coordinate system
$(t_1,t_2,\cdots,t_{2n-1},r)$ such that $A^{\theta} =
\frac{\partial}{\partial t_1}$, and $r$ is a defining function for
$\D$. $V$ can be assumed to be contained in this coordinate patch. We
will follow \cite{C78} to show that the integral curve of
$A^{\theta}$ from $t_1 = 0$ to $t_1 = t_0$, and hence $P$, can be
separated by the level set of a strictly plurisubharmonic function
from any compact subset of $\overline{\D} \setminus V$. This
contradiction (to $P \in \widehat{K} \setminus K$ for some compact
$K$) will complete the proof. As in \cite{C78}, consider the
auxiliary function  
$g(t_1,t_2,\cdots,t_{2n-1},r)=(t_1+\frac{1}{m})/(1+m^2(t_2^2+\cdots
+t_{2n-1}^2))$ 
and the sets
$S_c=\{(t_1,t_2,\cdots,t_{2n-1},0):g(t_1,t_2,\cdots,t_{2n-1},0)=c, 0
\leq t_1 \leq t_0 \}$. Choose $m$ sufficiently large such that $S_c
\subset \subset V\cap b\D $ for all $c$ with $1/m \leq c\leq
t_0+1/m$. By shrinking $V$ if necessary we may assume that any point
in the set $\{z\in b\D \cap V:t_1(z)\geq t_0\}$ is strongly
pseudoconvex. In \cite{C78}, $g$ is modified, but the modification is
specified to be of order $r^{2}$. Here, it is important to also have
a term that is of order $r$. Specifically, set $h=g+\mu r+\nu r^2$
for $\mu,\nu$ positive numbers to be determined. We have for $\tau >
0$ and $W
= \sum_{j=1}^{n}w_{j}\partial/\partial z_{j}$ (at boundary points,
where $r=0$)
$L_{e^{\tau h}}(W,\overline{W})=\tau e^{\tau
h}\Big( L_{g}(W,\overline{W})+\mu L_{r}(W,\overline{W}) +
\tau |W(g) + \mu W(r)|^{2}+ 2\nu |W(r)|^{2} \Big)$. 
Set $X_n = \big(1/\sum_{j=1}^{n}|\partial r/\partial
z_j|^{2}\big)\sum_{j=1}^{n}(\partial r/\partial 
\overline{z_j})\partial / \partial z_j)$ (so that $X_n(r)=1$). We
express $W$ in terms of 
the basis  $\{X_1, X_2, \cdots , X_{n-1}, X_n \}$ as 
$W=\sum_{j=1}^n\alpha_jX_j,$  
and apply the inequalities $2ab \leq \varepsilon
a^{2}+(1/\varepsilon)b^{2}$  
for $\varepsilon > 0$, $|a+b|^{2} \geq
\frac{1}{2}|a|^{2} - |b|^{2}$, and $\big|\sum_{j=1}^{n}a_j\big|^{2}
\leq n\sum_{j=1}^{n}|a_j|^{2}$, to obtain 
\begin{eqnarray*}
\frac{1}{\tau}e^{-\tau h}L_{e^{\tau h}}(W,\overline{W})&\geq
&\frac{\tau}{2}\left|X_1(g)\right|^{2}|\alpha_1|^{2} 
+ \sum_{j=2}^{n-1}\left(\mu L_{r}(X_j,\overline{X_j})-\tau
n\left|X_j(g)\right|^{2}\right)|\alpha_j|^{2}\\
&+& 2\nu |\alpha_n|^{2}-\left(\tau n\left|X_n(g)\right|^{2}+\tau
n\mu^{2}+
\mu\left|L_{r}(X_n,\overline{X_n})\right|\right)|\alpha_n|^{2}\\
&-& s.c.\sum_{j=1}^{n}|\alpha_j|^{2} - l.c.(\mu)|\alpha_n|^{2} + 
L_{g}(W,\overline{W}) \ ,
\end{eqnarray*}
where $s.c.$ denotes a small constant, and $l.c.(\mu)$ denotes a large
constant whose size depends on $\mu$. Note that on $\overline{V} \cap
b\D$, $|X_1(g)| \approx |A(g)|+|J(A)(g)| \geq
|A^{\theta}(g)|=|\partial g/\partial t_1|>0$ and
$L_{r}(X_j,\overline{X_j})>0$ for $2\leq j \leq n-1$. Therefore, we
can first choose $\tau$, then $\mu$, and then $\nu$, big enough so
that $L_{e^{\tau h}}(W,\overline{W}) \geq |W|^{2}$ on $\overline{V}
\cap b\D$. So $e^{\tau h}$ is strictly plurisubharmonic near
$\overline{V} \cap b\D$. From here on, the argument is exactly as in
\cite{C78}, p. 54-55; we only sketch it. Choose a smooth function
$\psi_s(t_1)$ which is identically $1$ for $t_1 \leq t_0$, and
identically $0$ for $t_1 \geq s$, where $s > t_0$ is chosen small
enough so that the level sets $\{\psi_{s}e^{\tau h} = c \}$ are
contained in $V$ for $1/m \leq c \leq t_0+1/m$. To deal with the
direction transverse to the boundary, Catlin applies his construction
in \cite{C78}, Theorem 3.1.6 (see also \cite{C80}, Proposition
3.1.6). In our situation, this construction yields a strictly
plurisubharmonic function on a neighborhood $V_1$ of $P$ with
$\{F_{A^{\theta}}^{t}(P): 0 \leq t \leq t_0 \} \subset\subset V_1
\subset\subset V$, whose superlevel set determined by $P$ is a
compact subset of $V_1$. Composition with a suitable convex
increasing function finally results in a plurisubharmonic function
defined on all of $\overline{\D}$ (by extension by $0$) that separates
$P$ from
any compact subset of $\overline{\D} \setminus V$. This completes the
proof of (the nontrivial direction of) Proposition \ref{hull}.

\bigskip
\smallskip

\bibliographystyle{amsplain}

\bigskip
\bigskip
\bigskip
\end{document}